\documentstyle[amssymb,12pt]{article}

\begin{document}

\author{Mircea Crasmareanu\thanks{%
Partially supported by the CEex Grant 05-D11-84}}
\title{Last multipliers for multivectors with applications to
Poisson geometry}
\date{}
\maketitle

\begin{abstract}
The theory of the last multipliers as solutions of the Liouville's
transport equation, previously developed for vector fields, is
extended here to general multivectors. Characterizations in terms
of Witten and Marsden differentials are reobtained as well as the
algebraic structure of the set of multivectors with a common last
multiplier, namely Gerstenhaber algebra. Applications to Poisson
bivectors are presented by obtaining that last multipliers count
for ''how far away'' is a Poisson structure from being exact with
respect to a given volume form. The notion of exact Poisson
cohomology for an unimodular Poisson structure on $I\!\!R^{n}$ is
introduced.
\end{abstract}



\vspace{.3cm}

\section*{Introduction}

In January 1838, Joseph Liouville(1809-1882) published a note, \cite{j:l},
on the time-dependence of the Jacobian of the "transformation" exerted by
the solution of an ODE on its initial condition. In modern language, if $%
A=A(x)$ is the vector field corresponding to the given ODE and
$m=m(t,x)$ is a smooth function (depending also on time $t$), then
the main equation of the cited paper is:
$$
\frac{dm}{dt}+m\cdot divA=0 \eqno(LE)
$$
called, by then, the {\it Liouville equation}. Some authors use the name
{\it ge\-ne\-ra\-lized Liouville equation}, \cite{g:e}, but we prefer to
name it the {\it Liouville equation of transport} (or {\it of continuity}).
This equation is a main tool in statistical mechanics where a solution is
called a {\it probability density function}, \cite{g:u}.

The notion of the {\it last multiplier}, introduced by Carl Gustav
Jacob Jacobi (1804-1851) around 1844, was treated in details in
{\it Vorlesugen \"{u}ber Dynamik}, edited by R. F. A. Clebsch in
Berlin in 1866. Thus, sometimes it has been used under the name of
{\it Jacobi multiplier}. Since then, this tool for understanding
ODE's was intensively studied by mathematicians in the usual
Euclidean space ${I\!\!R}^{n}$, as can be seen in the bibliography
of \cite{m:c1}, \cite{n:c1}-\cite{n:c4}. For all those interested
in historical aspects, an excellent survey can be found in
\cite{b:e}.

Several geometrical aspects of the last multipliers viewed as
autonomous, i.e. time-independent, solutions of LE are derived in
two papers by the same author: \cite{m:c1}, \cite{m:c2}. Our study
has been inspired by the results presented in \cite{o:z} using the
calculus on manifolds especially the Lie derivative, a well-known
tool for the geometry of vector fields.

The aim of the present paper is to extend this theory of the last
multipliers from vector fields to general multivectors by means of
the {\it curl operator}. This operator, a conjugate of usual
exterior derivative with respect to contraction of a given volume
form, was introduced by J.-L. Koszul in Poisson geometry
\cite{k:l} and is detailed in Chapter 2 of \cite{d:tz} and Section
2 of \cite{y:z}.

Since the Poisson multivectors are most frequently used, a Poisson
bracket is added to our study and we show that the last
multipliers are a measure of ''how far away'' is a Poisson
structure from being exact regarding the given volume form. Exact
Poisson structures are the theme of papers \cite{c:mm2} and
\cite{y:z} and form a remarkable class of Poisson structures
closed to symplectic structures as it is pointed out in \cite{a:w}
and the second paper cited above and proved in our Section 3.
There are other two important
features of these Poisson structures: \\
a) in \cite[p. 149]{d:tz} the problem of classification of quadratic
Poisson structures is reduced to the problem of classification of
exact quadratic Poisson structures and linear vector fields which preserve them,\\
b) \cite[Remark 3.2.]{y:z}: {\it in dimension 3 any Hamiltonian
vector field associated to an exact Poisson structure is
completely integrable}. \\
Let us remark that previously, in \cite{c:mm1}, the same notion was called {\it %
locally exact}.

The paper is structured as follows. The first section recalls the
definition of last multipliers and some previous results.
Characterizations in terms of other types of differentials than
the usual exterior derivative, namely Witten and Marsden, are
recalled from \cite{m:c2}. For a fixed smooth function $m$, the
set of vector fields admitting $m$ as last multiplier is shown to
be a Lie subalgebra of the Lie algebra of vector fields.

The next section is devoted to the announced extension to multivectors and
the previous results regarding Marsden and Witten differentials are
reobtained in this extended framework. Several consequences with respect to
the Schouten bracket on multivectors are derived including the extension of
final result from last paragraph.

In the following section the Poisson case is discussed and local expressions
for the main results of this section are provided in terms of the bivector $%
\pi $ defining the Poisson bracket. Again, last multipliers count
for the ''deformation'' from exactness of a given Poisson
structure. Two concrete examples (two-dimensional Poisson
structures and Lie-Poisson structures) are discussed and some
results of \cite{y:z} are reobtained in this way.

The last section is dedicated to a new notion namely {\it exact
Poisson co\-ho\-mo\-logy for an unimodular Poisson structure} in
$I\!\!R^{n}$. It is an open problem both the computation of this
cohomology and the relation with classical Poisson cohomology. For
this last theory details appear in \cite{d:tz} and \cite {va:i}.

{\bf Acknowledgments} The author expresses his thanks to ??? and
??? for several useful remarks.

\section{Last multipliers for vector fields}

Let $M$ be a real, smooth, $n$-dimensional manifold, $C^{\infty }\left(
M\right) $ the algebra of smooth real functions on $M$, ${\cal X}\left(
M\right) $ the Lie algebra of vector fields and $\Lambda ^{k}\left( M\right)
$ the $C^{\infty }\left( M\right) $-module of $k$-differential forms, $0\leq
k\leq n$. Assume that $M$ is orientable with the fixed volume form $V\in
\Lambda ^{n}\left( M\right) $.

Let:
\[
\stackrel{.}{x}^{i}\left( t\right) =A^{i}\left( x^{1}\left( t\right) ,\ldots
,x^{n}\left( t\right) \right) ,1\leq i\leq n
\]
be an ODE system on $M$ defined by the vector field $A\in {\cal
X}\left(
M\right) ,A=\left( A^{i}\right) _{1\leq i\leq n}$ and let us consider the $%
\left( n-1\right) $-form $\Omega_A =i_{A}V\in \Lambda ^{n-1}\left( M\right) $%
.

\medskip

{\bf Definition 1.1}(\cite[p. 107]{f:f}, \cite[p. 428]{o:z}) The function $%
m\in C^{\infty }\left( M\right) $ is called a {\it last multiplier} of the
ODE system generated by $A$, ({\it last multiplier }of $A$, for short) if $%
m\Omega_A $ is closed:
$$
d\left( m\Omega_A \right) :=\left( dm\right) \wedge \Omega_A +md\Omega_A =0.%
\eqno\left( 1.1\right)
$$

For example, in dimension $2$, the notions of the last multiplier
and integrating factor are identical and Sophus Lie suggested a
method to associate a last multiplier to every symmetry vector
field of $A$ (Theorem 1.1 in \cite[p. 752]{h:s}). Lie's method is
extended to any dimension in \cite {o:z}.

Characterizations of last multipliers can be obtained in terms of
Witten's differential \cite {w:e} and Marsden's differential
\cite[p. 220]{m:a}. If $f\in C^{\infty }\left( M\right) $ and
$t\geq 0$, Witten deformation of the usual differential
$d_{tf}:\Lambda ^{\ast }\left( M\right) \rightarrow \Lambda ^{\ast
+1}\left( M\right) $ is defined by:
\[
d_{tf}=e^{-tf}de^{tf}
\]
which means \cite{w:e}:
\[
d_{tf}\left( \omega \right) =tdf\wedge \omega +d\omega .
\]
Hence, $m$ is a last multiplier if and only if:
\[
d_{m}\Omega_A =\left( 1-m\right) d\Omega_A
\]
i.e. $\Omega_A $ belongs to the kernel of the differential operator $%
d_{m}+\left( m-1\right) d:\Lambda ^{n-1}\left( M\right) \rightarrow \Lambda
^{n}\left( M\right) $. Marsden differential is $d^{f}:\Lambda ^{\ast }\left(
M\right) \rightarrow \Lambda ^{\ast +1}\left( M\right) $ defined by:
\[
d^{f}\left( \omega \right) =\frac{1}{f}d\left( f\omega \right)
\]
and $m$ is a last multiplier if and only if $\Omega_A $ is
$d^{m}$-closed.

The following characterization of the last multipliers will be
useful:

\medskip

{\bf Lemma 1.2}(\cite[p. 428]{o:z}) (i) $m\in C^{\infty }\left( M\right) $
{\it is a last multiplier for} $A$ {\it if and only if}:
$$
A\left( m\right) +m\cdot div_{V}A=0\eqno\left( 1.2\right)
$$
{\it where $div_{V}A$ is the divergence of $A$ with respect to volume form $%
V $}.

(ii) {\it Let $0\neq h\in C^{\infty }\left( M\right) $ such that:
$$
L_{A}h:=A\left( h\right) =\left( div_{V}A\right) \cdot h\eqno\left(
1.3\right)
$$
Then $m=h^{-1}$ is a last multiplier for $A$}.

\medskip

{\bf Remarks} {\bf 1.3} (i) Equation $(1.2)$ is exactly the time-independent
version of LE from the Introduction. An important feature of equation $%
\left( 1.2\right) $ is that it does not always admit solutions
\cite[p. 269]{g:fc}. \\
(ii) In the terminology of \cite[p. 89]{b:e}, a function {\it h}
satisfying (1.3) is called an {\it inverse multiplier}. \\
(iii) A first result given by $\left( 1.2\right) $ is the
characterization of last multipliers for divergence-free vector
fields: $m\in C^{\infty }\left( M\right) $ {\it is a last
multiplier for the divergenceless vector field }$A$ {\it if and
only if} $m$ {\it is a first integral of} $A$. The importance of
this result is shown by the fact that three remarkable classes of
divergence-free vector fields are provided by: Killing vector
fields in Riemannian geometry, Hamiltonian vector fields in
symplectic geometry and Reeb vector fields in contact geometry.
Also, there are many equations of mathematical physics
corresponding to the vector
fields without divergence. \\
(iv) For the general case, namely $A$ is not divergenceless, there
is a strong connection between the first integrals and the last
multipliers as well. Namely, from properties of Lie derivative,
the ratio of two last multipliers is a first integral and
conversely, the product between a first integral and a last
multiplier is a last multiplier. So, denoting $FInt(A)$ the set of
first integrals of $A$, since $FInt(A)$ is a subalgebra in
$C^{\infty }(M)$ it results that the set of last multipliers for
$A$ is a $FInt(A)$-module. \\
(v) Recalling formula:
$$
div_{V}\left( fX\right) =X\left( f\right) +fdiv_{V}X\eqno(1.4)
$$
it follows that $m$ is a last multiplier for $A$ if and only if
the vector field $mA$ is with null divergence i.e. $div_{V}\left(
mA\right) =0$. Thus, the set of last multipliers is a ''measure of
how far away'' is $A$ from being divergence-free.

\medskip

An important structure generated by a last multiplier is given by:

\medskip

{\bf Proposition 1.4} {\it Let} $m\in C^{\infty }\left( M\right) $ {\it be
fixed. The set of vector fields admitting} $m$ {\it as last multiplier is a
Lie subalgebra in} ${\cal X}\left( M\right) $.

\medskip

{\bf Proof} Let $X$ and $Y$ be vector fields with the required property.
Since \cite[p. 123]{m:r}:
\[
div_{V}\left[ X,Y\right] =X\left( div_{V}Y\right) -Y\left( div_{V}X\right)
\]
one has:
\begin{eqnarray*}
\left[ X,Y\right] \left( m\right) +mdiv_{V}\left[ X,Y\right] &=&\left(
X\left( Y\left( m\right) \right) +mX\left( div_{V}Y\right) \right) - \\
-\left( Y\left( X\left( m\right) \right) +mY\left( div_{V}X\right) \right)
&=&\left( -div_{V}Y\cdot X\left( m\right) \right) -\left( -div_{V}X\cdot
Y\left( m\right) \right) = \\
&=&div_{V}Y\cdot mdiv_{V}X-div_{V}X\cdot mdiv_{V}Y=0. \\
&&\hfill \square
\end{eqnarray*}

\section{Last multipliers for multivectors}

Denote by ${\cal X}^{k}(M)$ the $C^{\infty }(M)$-module of $k$-vector fields, $%
1\leq k\leq n$ and fix $A\in {\cal X}^{k}(M)$. The multivector $A$ defines
the map $i_{A}:\Lambda^{p}(M)\rightarrow \Lambda^{p-k}(M)$ given by: \newline
$\cdot <i_A\omega, B>=<\omega, A\wedge B>$ for every $B\in {\cal X}^{p-k}(M)$
with $<,>$ the natural duality between forms and multivectors and $\wedge $
the Grassmann wedge product on $\bigoplus\limits_{k=1}^{n}{\cal X}^{k}\left(
M\right) $, if $p\geq k$, \newline
$\cdot \ i_A\omega =0$ if $p<k$.

It follows that on $(M, V)$ lives the map:
$$
V^\flat :{\cal X}^k(M)\rightarrow \Lambda^{n-k}(M), \qquad V^\flat(A)=i_AV, %
\eqno(2.1)
$$
which is a $C^{\infty }(M)$-isomorphism between ${\cal X}^k(M)$ and $%
\Lambda^{n-k}(M)$, for $0\leq k\leq n$. The inverse map of $V^\flat $ is
denoted $V^\natural :\Lambda^{n-k}(M)\rightarrow {\cal X}^k(M)$.

\medskip

{\bf Definition 2.1}(\cite[p. 70]{d:tz}) The map $D_{V}:{\cal X}%
^{k}(M)\rightarrow {\cal X}^{k-1}(M)$:
$$
D_{V}=V^{\natural }\circ d\circ V^{\flat },\eqno(2.2)
$$
is called {\it the curl operator} with respect to the volume form $V$. So,
if $A\in {\cal X}^{k}(M)$ then $D_{V}A$ is called {\it the curl of} $A$.

\medskip

{\bf Example 2.2}(\cite[p. 70]{d:tz}) If $k=1$ then $D_{V}=div_{V}$. Indeed,
if $A\in {\cal X}(M) $ then:
\[
\left( D_{V}A\right) V=V^{\flat }\circ D_{V}\left( A\right) =d\circ V^{\flat
}\left( A\right) =d\circ i_{A}\left( V\right) =L_{A}V=\left( div_{V}A\right)
V.
\]

\medskip

Inspired by this example and relation $\left( 1.4\right) $ we
introduce here the main notion of this paper:

\medskip

{\bf Definition 2.3} The function $m\in C^{\infty }\left( M\right) $ is
called a {\it last multiplier} of $A\in {\cal X}^{k}(M)$ if:
$$
D_{V}(mA)=0.\eqno\left( 2.3\right)
$$

\medskip

Since $V^{\natural }$ is a $C^{\infty }(M)$-isomorphism between $\Lambda
^{n-k}(M)$ and ${\cal X}^{k}(M)$ it results that $\left( 2.3\right) $ means $%
d\left( V^{\flat }\left( mA\right) \right) =0$ i.e.:
$$
d\left( mV^{\flat }\left( A\right) \right) =0\eqno\left( 2.4\right)
$$
which is the natural extension of condition $\left( 1.1\right) $
from Definition 1.1. With the same computation as in the previous
section we derive the following equivalent characterizations of
last multipliers for $A\in {\cal X}^{k}(M)$: \newline $\cdot $ in
terms of Witten differential: $V^{\flat }\left( A\right) =i_{A}V$
belongs to the kernel of the differential operator $d_{m}+\left(
m-1\right)
d:\Lambda ^{n-k}\left( M\right) \rightarrow \Lambda ^{n-k+1}\left( M\right) $%
, \newline
$\cdot $ in terms of Marsden differential: $V^{\flat }\left( A\right)
=i_{A}V $ is $d^{m}$-closed with $d^{m}:\Lambda ^{k}\left( M\right)
\rightarrow \Lambda ^{k+1}\left( M\right) $ as in Section 1.

\medskip

From the $C^{\infty }(M)$-linearity of $V^{\flat }$ we have $V^{\flat
}\left( mA\right) =mV^{\flat }\left( A\right) =\left( mV\right) ^{\flat
}\left( A\right) $ and then $\left( mV\right) ^{\natural }=\frac{1}{m}%
V^{\natural }$ (we suppose $m>0$ everywhere). It follows:
$$
mD_{mV}\left( A\right) =V^{\natural }\circ d\circ V^{\flat }\left( mA\right)
=D_{V}\left( mA\right) \eqno\left( 2.5\right)
$$
which yields:

\medskip

{\bf Proposition 2.4} $m\in C^{\infty }\left( M\right) $ {\it is a last
multiplier of} $A\in {\cal X}^{k}(M)$ {\it if and only if}:
$$
D_{mV}\left( A\right) =0.\eqno\left( 2.6\right)
$$

\medskip

The last formula has some important consequences, all in terms of
an operation on $\bigoplus\limits_{k=1}^{n}{\cal X}^{k}\left(
M\right) $ called {\it Schouten bracket} $\left[ ,\right] $ which
is a natural ge\-ne\-ra\-lization of Lie bracket from ${\cal
X}\left( M\right) $ and generates a {\it Gersternhaber algebra}
structure on the set of multivectors, \cite{y:ks}. For details
regarding this bracket see \cite{d:tz}, \cite{va:i}. The first
corollary of $\left( 2.6\right) $ is a formula for the curl:

\medskip

{\bf Proposition 2.5} {\it If} $m\in C^{\infty }\left( M\right) $
{\it is a non-vanishing last multiplier of} $A\in {\cal X}^{k}(M)$
{\it then the curl of} $A$ {\it can be expressed in terms of the
Schouten bracket}:
$$
D_{V}A=-\left[ A,\ln |m|\right] .\eqno\left( 2.7\right)
$$

\smallskip

{\bf Proof} Is a direct consequence of formula $\left( 2.90\right) $ from
\cite[p. 71]{d:tz}:
$$
D_{mV}A=D_{V}A+\left[ A,\ln |m|\right] . \eqno \Box
$$

\medskip

A second formula relates the Schouten bracket with the product $\wedge $ of $%
\bigoplus\limits_{k=1}^{n}{\cal X}^{k}\left( M\right) $. After \cite[Th.
2.6.7 p. 71]{d:tz} if $A$ is an $a$-multivector and $B$ is a $b$-multivector
then:
$$
\left[ A,B\right]=\left( -1\right) ^{b}D_{V}\left( A\wedge B\right) -\left(
D_{V}A\right) \wedge B-\left( -1\right) ^{b}A\wedge \left( D_{V}B\right) .%
\eqno(2.8)
$$

\smallskip

{\bf Corollary 2.6} {\it Let} $m\in C^{\infty }\left( M\right) $ {\it be a
last multiplier for both} $A$ {\it and} $B$. {\it Then} $m$ {\it is a last
multiplier for} $A\wedge B$ {\it if and only if} $A$ {\it and} $B$ {\it %
Schouten-commutes i.e. their Schouten bracket vanishes}: $\left[ A,B\right]
=0$.

\medskip

Another consequence of $\left( 2.6\right) $ is a straightforward
generalization of Proposition 1.4:

\medskip

{\bf Theorem 2.7} {\it Let} $m\in C^{\infty }\left( M\right) $
{\it be fixed. The set of multivectors admitting} $m$ {\it as last
multiplier is a Gerstenhaber subalgebra in}
$\bigoplus\limits_{k=1}^{n}{\cal X}^{k}\left( M\right) $.

\medskip

{\bf Proof} The curl operator is, up to a sign, a derivation of the Schouten
bracket, namely \cite[p. 71]{d:tz}:
$$
D_{V}\left[ A,B\right]=\left[ A,D_{V}B\right] +\left( -1\right) ^{b-1}\left[
D_{V}A,B\right] .\eqno(2.9)
$$
This relation combined with $(2.6)$ gives the conclusion. $\Box $

\medskip

{\bf Definition 2.8}(\cite{y:z}) The multivector $A$ is called {\it exact}
with respect to the volume form $V$ if $D_{V}\left( A\right) =0$.

\medskip

{\bf Remarks 2.9} (i) It follows from $\left( 2.3\right) $ that the set of
last multipliers of $A$ is a ''measure of how far away'' is $A$ from being
exact. \newline
(ii) Equation $(2.8)$ gives that if $A$ and $B$ are exact multivectors then $%
A\wedge B$ is exact if and only if they Schouten-commutes.
\newline (iii) Using again $(2.9)$ it results that the set of
exact multivectors is a Schouten subalgebra in
$\bigoplus\limits_{k=1}^{n}{\cal X}^{k}\left( M\right) $.

\medskip

{\bf Example 2.10} From \cite{m:t} the volume form $V$ yields a
Nambu multivector, \cite[p. 160]{d:tz}, $A_{V}\in {\cal
X}^{n}\left( M\right) $; if $\left( x^{1},\ldots ,x^{n}\right) $
is a local chart on $M$ such that $V=fdx^{1}\wedge \ldots \wedge
dx^{n}$ then $A_{V}=\frac{1}{f}\frac{\partial }{\partial
x^{1}}\wedge \ldots \wedge \frac{\partial }{\partial x^{n}}$. A
straightforward computation gives that $A_{V}$ is exact with
respect to $V$.

\medskip

{\bf Remark 2.11} Let $f\in C^{\infty}(M)$ and $A$ an
$a$-multivector. From $(2.8)$ and $D_V(f)=0$ we get:
$$
\left[A,f\right]=D_V(fA)-fD_V(A)
$$
and then $D_V(fA)=fD_V(A)$ if and only if $f$ is a {\it Casimir}
of $A$ i.e. $\left[A,f\right]=0$. Connecting this with Remarks
2.9.\ (ii) we derive:

\medskip

{\bf Proposition 2.12} {\it If $A$ is exact then $fA$ is exact if
and only if $f$ is a Casimir function of} $A$.

\section{Last multipliers for Poisson bivectors}

Let us assume that $M$ is endowed with a Poisson bracket $\{,\}$ induced by
the Poisson bivector $\pi \in {\cal X}^{2}\left( M\right) $. Let $f\in
C^{\infty }\left( M\right) $ and $A_{f}\in {}\left( M\right) $ be the
associated {\it Hamiltonian vector field} of the {\it Hamiltonian} $f$, \cite
{m:r}.

Given the volume form $V$ there exists a unique vector field $X_{\pi ,V}$,
called the {\it modular vector field}, so that \cite{k:l}, \cite{a:w}:
$$
div_{V}A_{f}=X_{\pi ,V}\left( f\right) .\eqno\left( 3.1\right)
$$

From Proposition 1 of \cite[ p. 4]{d:f} we have:
$$
X_{\pi ,V}=D_{V}\left( \pi \right) .\eqno(3.2)
$$

\smallskip

{\bf Definition 3.1} The triple $\left( M,\pi ,V\right) $ is called \cite
{a:w} {\it unimodular} if $X_{\pi ,V}$ is a Hamiltonian vector field, $%
A_{\rho }$ of $\rho \in C^{\infty }\left( M\right) $. The triple $\left(
M,\pi ,V\right) $ is called \cite{c:mm2}, \cite{y:z} {\it exact} if $X_{\pi
,V}$ is identically zero.

\medskip

Let us introduce:

\medskip

{\bf Definition 3.2} The function $m\in C^{\infty }\left( M\right) $ is
called a {\it last multiplier} of $\left( M,\pi ,V\right) $ if:
$$
D_{V}\left( m\pi \right) =0\eqno(3.3)
$$
equivalently:
$$
D_{mV}\left( \pi \right) =0.\eqno\left( 3.4\right)
$$

\medskip

It results that the set of the last multipliers of $\left( M,\pi
,V\right) $ is a ''measure of how far away'' is $\left( M,\pi
,V\right) $ from being exact and the characterization:

\medskip

{\bf Proposition 3.3} $m\in C^{\infty }\left( M\right) $ {\it is a last
multiplier} of $\left( M,\pi ,V\right) $ {\it if and only if}:
$$
X_{\pi, mV}=0.\eqno(3.5)
$$

\medskip

{\bf Example 3.4} i) Poisson structures induced by symplectic
structures are exact. This statement appears in the introduction
of \cite{y:z} and we provide here a proof using \cite{a:w}(or item
1 of Remark 2.3. from \cite {y:z}): a Poisson structure is exact
with respect to $V$ if and only if $V$ is invariant of any
Hamiltonian vector field $A_{f}$. But in symplectic geometry this
is a well-known fact. \\
ii) A condition for a quadratic Poisson structure on $I\!\!R^{3}$
to be exact is given in Example 5.6.8. from \cite[p. 149]{d:tz}.

\medskip

The two notions of Definition 3.1 are equivalent as it is pointed
out in \cite {c:mm2}. Moreover, in the MR review of \cite{y:z} it
is put in evidence that at local level there is no problem about
the dependence of volume form $V$. So, in the following we work in
local coordinates. Let $\left( x^{1},\ldots ,x^{n}\right) $ be a
local chart on $M$ such that $V=dx^{1}\wedge \ldots \wedge dx^{n}$
and the bivector $\pi $ of $\left( M,\{,\}\right) $ is: $\pi
=\sum\limits_{i<j}\pi ^{ij}\frac{\partial }{\partial x^{i}}\wedge \frac{%
\partial }{\partial x^{j}}$. Denoting $\pi ^{i}=\sum\limits_{j=1}^{n}\frac{%
\partial \pi ^{ij}}{\partial x^{j}}$ we have \cite[Proposition 1, p. 4]{d:f}%
, \cite{c:mm2}:
$$
X_{\pi, V}=\sum\limits_{i=1}^{n}\pi ^{i}\frac{\partial }{\partial x^{i}}\eqno%
\left( 3.6\right)
$$
and then, Proposition 3.3 becomes:

\medskip

{\bf Proposition 3.5} $m\in C^{\infty }\left( M\right) $ {\it is a last
multiplier for }$\left( M,\pi ,V\right) ${\it \ if and only if}:
$$
\pi _{m}^{i}:=\sum\limits_{j=1}^{n}\frac{\partial \left( m\pi ^{ij}\right) }{%
\partial x^{j}}=0,\qquad 1\leq i\leq n.\eqno(3.7)
$$

\smallskip

{\bf Examples 3.6:}

{\bf 3.6.1}

After \cite[p. 31]{va:i} the bivector $\pi =h\left( x,y\right) \frac{%
\partial }{\partial x}\wedge \frac{\partial }{\partial y}$ defines a Poisson
structure on ${I\!\!R}^{2}$. So, $\pi ^{12}=-\pi ^{21}=h$ and then $\left(
3.7\right) $ becomes:
\[
\frac{\partial \left( mh\right) }{\partial y}=-\frac{\partial \left(
mh\right) }{\partial x}=0
\]
with the obvious solution $m_{\pi }=\frac{C}{h}$ (if we suppose $h>0$
everywhere), where $C$ is a real constant. Therefore, on the Poisson manifold $%
\left( {I\!\!R}^{2},\pi \right) $ above, the function $C/h$ is a
last multiplier.

In this way we reobtain part (a) of Theorem 3.2. from \cite{y:z}
that any smooth 2-dimensional Poisson structure is exact if and
only if it is constant; indeed the exact Poisson
$m_{\pi }\cdot \pi =C\frac{\partial }{\partial x}\wedge \frac{\partial }{%
\partial y}$ is constant. Also, the second phrase of Remark 3.2. item 3):
{\it the set of exact 2-dimensional Poisson structures is a
1-dimensional space isomorphic with} ${I\!\!R}$ is also verified.

\medskip

{\bf 3.6.2 Lie-Poisson structures}

\medskip

The interest for this example is pointed out in \cite{y:z}:\ {\it %
Lie-Poisson structures play important roles in studying normal forms for a
class of Poisson structures.}

Let ${\cal G}$ be an $n$-dimensional Lie algebra with a fixed basis $%
B=\{e_{i}\}_{1\leq i\leq n}$ and let $B^{\ast }=\{e^{i}\}$ be the dual basis
on the dual ${\cal G}^{\ast }$. Recall the definition of {\it structure
constants} of ${\cal G}$:
\[
\left[ e_{i},e_{j}\right] =c_{ij}^{k}e_{k}.
\]
Then, on ${\cal G}^{\ast }$ we have the so-called {\it Lie-Poisson structure}
given by \cite[p. 31]{va:i}:
$$
\pi ^{ij}\left( x_{u}e^{u}\right) =c_{ij}^{k}x_{k}.\eqno\left( 3.8\right)
$$
We get:
$$
\pi _{m}^{i}=\sum\limits_{j=1}^{n}c_{ij}^{k}\frac{\partial \left(
mx_{k}\right) }{\partial x_{j}}\eqno\left( 3.9\right)
$$

{\bf Particular case: n=2}

Although from the previous example we know all about the
2-dimensional case it is interesting to reobtain the conclusion
within this example. The structure relations
$\left[ e_{1},e_{1}\right] =[e_{2},e_{2}]=0,\left[ e_{1},e_{2}%
\right] =c_{12}^{1}e_{1}+c_{12}^{2}e_{2}$ yield:
$$
\left\{
\begin{array}{c}
\pi _{m}^{1}=c_{12}^{2}m+\frac{\partial m}{\partial y}\left(
c_{12}^{1}x+c_{12}^{2}y\right) \\
\pi _{m}^{2}=-c_{12}^{1}m-\frac{\partial m}{\partial x}\left(
c_{12}^{1}x+c_{12}^{2}y\right)
\end{array}
\right. .\eqno\left( 3.10\right)
$$

Supposing ${\cal G}$\ nontrivial (i.e. $\left( c_{12}^{1}\right)
^{2}+\left( c_{12}^{2}\right) ^{2}>0$) there result three cases:
\newline I) $c_{12}^{1}\cdot c_{12}^{2}\neq 0$ i.e.
$h=c_{12}^{1}x+c_{12}^{2}y$. From the system $(3.7)$ $\pi
_{m}^{1}=\pi _{m}^{2}=0$ we have:
$$
c_{12}^{2}\frac{\partial m}{\partial x}-c_{12}^{1}\frac{\partial m}{\partial
y}=0\eqno\left( 3.11\right)
$$
with solution $m=A\left( \frac{x}{c_{12}^{2}}+\frac{y}{c_{12}^{1}}\right) +B$
which replaced in $\left( 3.10\right) $ yields $A=B=0$. In conclusion, the
last multiplier of $\pi $ for this case is zero and the associated Poisson
structure is trivial (hence exact). \newline
II) $c_{12}^{2}=0$ (i.e. $h=c_{12}^{1}x$) with solution $m=m\left( x\right) $
of $\left( 3.11\right) $. Inserting this function in $\left( 3.10_{2}\right)
$\ we get $m+x\cdot m^{\prime }=0$ with solution $m_{\pi }=\frac{C}{x}$.
\newline
III) $c_{12}^{1}=0$ (i.e. $h=c_{12}^{2}y$) with solution $m=m\left( y\right)
$ of $\left( 3.12\right) $. With the same computations as above it results $%
m_{\pi }=\frac{C}{y}$.

\section{Exact Poisson cohomology of unimodular Poisson structures}

Returning to the general case of Poisson structures in $I\!\!R^n$
let us point out an interesting consequence of $(2.8)$ and $(2.9)$
respectively:

\medskip

{\bf Proposition 3.7} i) {\it Let $X, Y\in {\cal X}(I\!\!R^n)$ be such that:
\newline
a) their wedge product $\pi =X\wedge Y$ is a Poisson structure,
\newline b) they Lie-commutes: $[X,Y]=0$. \newline c) they are
divergence-free. \newline Then $\pi $ is an unimodular Poisson
bivector}. \newline
ii) {\it Let $\pi $ be a Poisson structure and
$X\in {\cal X}(I\!\!R^n)$ such that their Schouten bracket $\left[
\pi ,X\right] $ is again a Poisson structure. If $\pi $ is
unimodular and $X$ is divergence-free then $\left[ \pi ,X\right]$
is unimodular}. \newline
iii) {\it Let $\pi $ be an unimodular
Poisson structure and $A$ an exact multivector. Then their
Schouten bracket $[\pi,A]$ is an exact multivector}.

\medskip

In the following suppose $\left( I\!\!R^{n},\pi \right) $ is an unimodular
Poisson manifold. Let us consider, after \cite[p. 39]{d:tz}, the map $\delta
_{\pi }:\bigoplus\limits_{k=1}^{n}{\cal X}^{k}\left( I\!\!R^{n}\right)
\rightarrow \bigoplus\limits_{k=1}^{n}{\cal X}^{k}\left( I\!\!R^{n}\right)
,\delta _{\pi }\left( A\right) =\left[ \pi ,A\right] $ (for a local expression
see \cite[Formula (4.8), p. 43]{va:i}) and let us denote $%
{\cal X}_{e}^{k}\left( I\!\!R^{n}\right) $ the set of exact $k$%
-multivectors. From the last item of the previous result and the fact that $%
\left( \bigoplus\limits_{k=1}^{n}{\cal X}^{k}\left( I\!\!R^{n}\right)
,\delta _{\pi }\right) $ is a complex \cite[p. 39]{d:tz}, it results a new
differential complex:
$$
\ldots \rightarrow {\cal X}_{e}^{k-1}\left( I\!\!R^{n}\right) \stackrel{%
\delta _{\pi }}{\rightarrow }{\cal X}_{e}^{k}\left( I\!\!R^{n}\right)
\stackrel{\delta _{\pi }}{\rightarrow }{\cal X}_{e}^{k+1}\left(
I\!\!R^{n}\right) \rightarrow \ldots \eqno\left( 4.1\right)
$$
which will be called {\it the exact Lichnerowicz complex}. Let us
call the cohomology of this complex {\it exact Poisson
cohomology}. Obviously, the exact Poisson cohomology is included
in the usual Poisson cohomology treated in detail in \cite{d:tz}
and \cite{va:i}.

Therefore we set the exact Poisson groups:
$$
H_{e}^{k}\left( I\!\!R^{n},\pi \right) =\frac{\ker \{\delta _{\pi
}:{\cal X}_{e}^{k}\left( I\!\!R^{n}\right) \rightarrow {\cal
X}_{e}^{k+1}\left( I\!\!R^{n}\right) \}}{Im\{\delta _{\pi }:{\cal
X}_{e}^{k-1}\left( I\!\!R^{n}\right) \rightarrow {\cal
X}_{e}^{k}\left( I\!\!R^{n}\right) \}} . \eqno\left( 4.2\right)
$$
$H_{e}^{k}\left( I\!\!R^{n},\pi \right) $ is a subgroup of the
group $H^{k}\left( I\!\!R^{n},\pi \right) $ of Poisson cohomology.
For example $H_{e}^{0}\left( I\!\!R^{n},\pi \right)=H^{0}\left(
I\!\!R^{n},\pi \right)$ which is the group of Casimir functions of
$\pi $, \cite[p. 40]{d:tz}.

\section*{Conclusions}

0) The last multipliers constitute a measure to count the
"perturbation" from exactness. So, this notion can be thought in
the framework of \cite{b:m}. \\
1) The theory of the last multipliers can be extended from vector
fields to general multivectors preserving a series of remarkable
characterizations and results.\\
2) An important structure generated by a last multiplier is of
algebraic nature: the set of multivectors with a prescribed last
multiplier is a Gerstenhaber subalgebra.\\
3) From the two previous remarks it results that a natural
extension of our theory seems to work on Lie algebroids using the
tools of \cite{g:mm} and \cite{y:ks}. Hence, a sequel paper
\cite{c:h} is forthcoming.

\vspace{.2cm}

\noindent Faculty of Mathematics \newline
University "Al. I. Cuza"\newline
Ia\c si, 700506\newline
Rom\^ania\newline
e-mail: mcrasm@uaic.ro \newline
\newline
\noindent http://www.math.uaic.ro/$\sim$mcrasm

\end{document}